\documentclass[12pt]{article}
\usepackage{amsmath,mathrsfs,amssymb,latexsym,amsthm,times}

\def\dd#1{{\,\rm d}#1}

\newtheorem{thm}{Theorem}

\newtheorem{prop}{Proposition}

\newcommand\samefootnotemark{\addtocounter{footnote}{-1}\footnotemark}

\title{Delange's Tauberian theorem and asymptotic normality of random
ordered factorizations of integers}
\author{{\sc Hsien-Kuei Hwang\thanks{Part of this work
       was done while both authors were visiting Institut
       Mittag-Leffler, Djursholm, Sweden. They
       thank the Institute for hospitality and support.}}\\
   Institute of Statistical Science\\
   Academia Sinica \\
   Taipei 115\\
   Taiwan \and
   {\sc Svante Janson}\samefootnotemark\\
   Department of Mathematics \\
   Uppsala University \\
   PO Box 480, SE-751 06\\
   Uppsala\\ Sweden
}
\date{February 19, 2009}

\begin{document}
\maketitle

\begin{abstract}

By a suitable shifting-the-mean parametrization at the Dirichlet
series level and Delange's Tauberian theorems, we show that the
number of factors in random ordered factorizations of integers is
asymptotically normally distributed.

\end{abstract}

\section{Introduction}

Let $\mathscr{P}$ be a fixed subset of $\{2,3,\dots\}$. Let $a(n)$
denote the number of different ways of writing $n$ as the product of
\emph{ordered} sequences $(n_1,\dots,n_k)$ of integers in
$\mathscr{P}$. Define $a(1)=1$. Let $A(N) := \sum_{1\le n\le N}
a(n)$. Assume that all $A(N)$ factorizations of an integer $\le N$
are equally likely; denote by $Y_N$ the random variable counting the
number of factors in a random factorization. We prove in this paper
that the distribution of $Y_N$ is asymptotically normal under very
general conditions on $\mathscr{P}$.

Denote by $\mathcal{P}(s)$ the Dirichlet series of $\mathscr{P}$
\[
    \mathcal{P}(s) := \sum_{n\in\mathscr{P}} n^{-s}.
\]
Assume the abscissa of convergence of $\mathcal{P}(s)$ is $\kappa$.
Then $\kappa=-\infty$ if $|\mathscr{P}|<\infty$ and $0\le \kappa\le
1$ if $|\mathscr{P}|=\infty$. Note that
$\mathcal{P}(\kappa)\le\infty$.

Our main result is as follows.

\begin{thm} \label{T:thm} Assume
$|\mathscr{P}|\ge2$ and\/ $1<\mathcal{P}(\kappa)\le \infty$; thus
there exists $\rho>\max\{\kappa,0\}$ such that
$\mathcal{P}(\rho)=1$. Let
\[
    \mu := -\frac1{\mathcal{P}'(\rho)},\quad\text{ and }\quad
    \sigma^2 := \mu^3\mathcal{P}''(\rho)-\mu.
\]
Then
\[
    \frac{Y_N-\mu\log N}{\sigma\sqrt{\log N}}
    \stackrel{d}{\to} \mathscr{N}(0,1),
\]
where $\stackrel{d}{\to}$ stands for convergence in distribution and
$\mathscr{N}(0,1)$ denotes the standard normal distribution. The
mean and variance of\/ $Y_N$ are asymptotic to $\mu\log N$ and
$\sigma^2\log N$, respectively.
\end{thm}
By Cauchy-Schwarz inequality, $\sigma^2>0$. We will indeed prove
convergence of all moments.

The case when $|\mathscr{P}|=1$, say $\mathscr{P}=\{d\}$, $d\ge2$,
is exceptional. In this case, $a(n)=1$ when $n$ is a power of $d$,
and $a(n)=0$ otherwise. Then $Y_N$ is uniformly distributed on the
integers $\{0, \dots,\lfloor \log_d N\rfloor\}$, and thus $Y_N/\log
N$ converges in distribution to a uniform distribution on $[0,1/\log
d]$, and therefore, $Y_N$ is not asymptotically normal; further, the
moment asymptotics differ from those in the theorem.

Ordered factorization problems in connection with that studied in
this paper have a long history, tracing back to at least MacMahon's
work (see \cite{MacMahon91}) in the early 1890's; later they were in
most publications referred to as Kalm\'ar's problem of
``factorisatio numerorum" (see \cite{Kalmar31,Hwang94}). Diverse
properties of such factorizations have then been widely
investigated, often in quite different contexts, one reason being
that ordered factorizations are naturally encountered in many
enumeration problems. For example, when $\mathscr{P}=\{2,3,\dots\}$,
$A(N)+1$ equals the permanent of the Redheffer matrix; see
\cite{Wilf04}; also they appeared as the lower bound of certain
biological problems; see \cite{NN93}. See also \cite{DHN08,KM05} for
more information and references.

The first paper dealing with general ordered factorizations beyond
the subset $\mathscr{P}=\{2,3,\dots\}$ similar to our setting was
Erd\H os \cite{Erdos41}, extending previous results by Hille
\cite{Hille36}; see also \cite{KKW93}. Asymptotic normality of the
special case of Theorem~\ref{T:thm} when $\mathscr{P}=\{2,3,\dots\}$
was treated in \cite{Hwang00}. In this case, $\rho\approx 1.7286$
being the unique root $>1$ that solves $\zeta(s)=2$, where $\zeta$
denotes Riemann's zeta function. The proof given there relies on the
determination of a zero-free region of the function
$1-z(\zeta(s)-1)$, which in turn involves deep estimates from
trigonometric sums (see also \cite{Lau01}). Such refined estimates
are for general $\mathscr{P}$ hard to establish. We replace this
estimate by applying purely Tauberian arguments of Delange (see
\cite{Delange54}), which require only analytic information of the
involved Dirichlet series on its half-plane of convergence.

We will use Dirichlet series and the method of moments and derive
asymptotic estimates for central moments of integral orders, which,
by Frechet-Shohat's moment convergence theorem, will suffice to
prove the theorem.
\begin{prop} For $k\ge0$ \label{T:pp1}
\begin{align} \label{T:cm}
    \lim_{N\to\infty}\mathbb{E}\left(\frac{Y_N-\mu\log N}
    {\sigma\sqrt{\log N}}\right)^k=
    \left\{\begin{array}{ll}
        \displaystyle \frac{k!}{(k/2)!2^{k/2}},
        &\text{if $k$ is even}\\
        0,& \text{if $k$ is odd},
    \end{array}\right.
\end{align}
\end{prop}

A straightforward application of the Tauberian theorem does not
provide precise asymptotics for central moments beyond the first due
to cancellation of major dominant terms and due to the fact that no
error term is generally available via application of Tauberian
arguments. The new idea we introduce in this paper is to take into
account the feature that the mean is logarithmic and to
shift-the-mean on the associated Dirichlet series, which nicely
incorporates the cancellations of higher central moments in a
surprisingly neat way. It thus avoids completely the messy
calculations and cancellations that the usual method of moments faces
when dealing with higher central moments. A similar idea was
previously applied to characterize the phase change of random
$m$-ary search trees, where a nonlinear differential equation with
an additive nature was encountered; see \cite{CH01}. However, the
tools used there are complex-analytic, in contrast to the purely
Tauberian ones used here.

In the special case when $\mathscr{P}= \{2,3,\dots\}$, our result
may be interpreted as saying that the property of the zero-free
region for the Dirichlet series $1-z(\zeta(s)-1)$ lies much deeper
than the asymptotic normality of the random variables $Y_N$.

In addition to the number-theoretic interest \emph{per se} of our
results, we believe that the approach we use here also offers
important methodological value for the study of similar problems. In
particular, not only is the use of the method of moments very
simple, but no analytic properties of the Dirichlet series beyond
the abscissa of convergence are needed, which largely simplifies the
analysis in many situations. For example, our approach can be
readily applied to other factorizations such as branching or cyclic
factorizations with algebraic or logarithmic singularities (see
\cite{Hwang94}). It can also be extended, coupling with suitable
Tauberian theorems, to deal with ordered factorizations in additive
arithmetic semigroups (see \cite{KKW93}) and components counts in
ordered combinatorial structures (see \cite{FS93}). Another possible
extension is to the analysis of Euclidean algorithms; see
\cite{BV05,FV98,Vallee97}.

In the periodic case when $\mathscr{P}\subseteq\{ d^k\,:\,k\ge1\}$
for some $d\ge2$, $\mathcal{P}(s)$ has period $2\pi i/\log d$ and
the usual Tauberian theorem does not apply. Instead we write
$\mathcal{P}(s) = \tilde{\mathcal{P}}(d^{-s})$, where
$\tilde{\mathcal{P}}(z) := \sum_{d^k\in\mathscr{P}} z^k$, and
transform the multiplicative nature of the problem into an additive
one on random compositions by taking logarithms. We replace the
Tauberian theorem by the singularity analysis of Flajolet and
Odlyzko (see \cite{FO90} or \cite[Chapter VI]{FS09}) in the proof;
the details are similar to the proof below, but simpler, so we omit
them. In the rest of the paper we thus assume, for every $d\ge2$,
\begin{align}\label{T:c1}
    \mathscr{P} \not\subseteq \{d^k\}_{k\ge1}.
\end{align}

\section{Dirichlet series, Delange's Tauberian theorem, and proofs}

\paragraph{Generating functions.} Let $a_m(n)$ denote the number of
ordered factorizations of $n$ into exactly $m$ factors. Then
\[
    \mathcal{P}(s)^m = \sum_{n\ge1} a_m(n) n^{-s},
\]
in formal power series sense; analytically, we can take $s$
satisfying $\mathcal{P}(\Re(s))<\infty$. Thus if
$e^{\Re(z)}\mathcal{P}(\Re(s))<1$, then by absolute convergence
\begin{align}\label{T:bgf}
    \sum_{n\ge1}n^{-s}\sum_{m\ge0}a_m(n) e^{mz} =
    \sum_{m\ge0} e^{mz} \mathcal{P}(s)^m
    =\frac{1}{1-e^z\mathcal{P}(s)}.
\end{align}

\paragraph{Delange's Tauberian theorem.} We need the following form
of Delange's Tauberian theorem (see \cite{Delange54} or \cite[Ch.\
III, Sec.\ 3]{Narkiewicz83}).
\begin{quote}
{\sl Let $F(s):=\sum_{n\ge1}\alpha(n) n^{-s}$ be a Dirichlet series
with nonnegative coefficients and convergent for $\Re(s)>\varrho>0$.
Assume (i) $F(s)$ is analytic for all points on $\Re(s)=\varrho$
except at $s=\varrho$; (ii) for $s\sim\varrho$, $\Re(s)>\varrho$,
\[
    F(s) = \frac{G(s)}{(s-\varrho)^{\beta}}+H(s)\qquad(\beta>0),
\]
where $G$ and $H$ are analytic at $s=\varrho$ with
$G(\varrho)\neq0$. Then
\begin{align}
    \sum_{n\le N} \alpha(n) \sim \frac{G(\varrho)}
    {\varrho \Gamma(\beta)} N^\varrho(\log N)^{\beta-1}.
    \label{T:DTT}
\end{align}
}
\end{quote}

\paragraph{Asymptotics of $A(N)$.} Taking $z=0$ in (\ref{T:bgf}), we
obtain the Dirichlet series for $a(n)=\sum_{m\ge0} a_m(n)$
\[
    \mathcal{A}(s) = \sum_{n\ge1} a_n n^{-s} =
    \frac1{1-\mathcal{P}(s)},
\]
as long as $\Re(s)>\rho$. Note that the non-periodicity assumption
(\ref{T:c1}) implies that $\mathcal{P}(s)\not=1$ for all $s$ with
$\Re(s)=\rho$ but $s\not=\rho$. Hence $\mathcal{A}(s)$ is not only
analytic in the open half-plane $\{s\, :\,\Re(s)>\rho\}$ but also on
the boundary $\{s\,:\,\Re(s)=\rho\}$ except at $s=\rho$. The same
holds true for all Dirichlet series we consider below.

Now for our $\mathcal{P}(s)$, since $\mathcal{P}'(\rho)
=-\sum_{n\in\mathscr{P}} n^{-\rho} \log n < 0$, we see that
$\mathcal{P}(s)$ has a simple zero at $s=\rho$, and thus
$\mathcal{A}(s)$ has a simple pole at $\rho$ with
\[
    \mathcal{A}(s)=\frac1{1-\mathcal{P}(s)}
    \sim \frac{-1}{\mathcal{P}'(\rho)(s-\rho)},
\]
as $s\to \rho$. Hence Delange's Tauberian theorem applies and we
obtain
\begin{align}\label{T:AN}
    A(N) =\sum_{n\le N} a(n) \sim R N^{\rho},\qquad
    R := -\frac1{\rho \mathcal{P}'(\rho)} = \frac{\mu}{\rho}.
\end{align}
Furthermore, we also have
\[
    \sum_{n\le N} a(n) (\log n)^k \sim
    \frac{\mu}{\rho} \,N^\rho (\log N)^k
    \sim A(N)(\log N)^k,
\]
either by repeating the same procedure for the Dirichlet series
\begin{equation}
  \label{aks}
    (-1)^k \mathcal{A}^{(k)}(s) = \sum_{n\ge1}a(n)
    (\log n)^k n^{-s} = (-1)^k\frac{\dd{}^k}
    {\dd{s}^k}\frac1{1-\mathcal{P}(s)},
\end{equation}
or by using directly (\ref{T:AN}). The estimate will be used later.

\paragraph{The expected value of $Y_N$.} By taking the derivative with
respect to $z$ on both sides of (\ref{T:bgf}), we obtain
\[
    \sum_{n\ge1} n^{-s} \sum_{m\ge0} ma_m(n)
    = \frac{\mathcal{P}(s)}{(1-\mathcal{P}(s))^2}.
\]
Delange's Tauberian conditions being easily checked as above, we
then obtain
\[
    \mathbb{E}(Y_N) = \frac1{A(N)}
    \sum_{n\le N}\sum_{m\ge0} ma_m(n)
    \sim \mu \log N.
\]

\paragraph{Shifting-the-mean at the Dirichlet series level.} For
higher central moments, the idea we will use can \emph{formally} be
described by using Perron's integral representation as follows
(using \eqref{T:bgf} and for simplicity assuming temporarily that
$N$ is not an integer).
\begin{align*}
    \mathbb{E}\left(e^{(Y_N-\mu\log N)z}\right) &=
    \frac{1}{2\pi iA(N)} \int_{c-i\infty}^{c+i\infty}
    \frac{N^{s-\mu z}}{s}\frac1{1-e^z\mathcal{P}(s)}\dd s \\
    &= \frac{1}{2\pi iA(N)} \int_{c-i\infty}^{c+i\infty}
    \frac{N^s}{s}\frac1{(1+\mu z/s)(1-e^z\mathcal{P}(s+\mu z))}\dd s,
\end{align*}
where $c$ is suitably chosen; the fact that the mean being of order
$\log N$ is crucial here. We then \emph{formally} expect that
\begin{align}\label{T:fsz}
    \mathbb{E}\left(Y_N-\mu\log N\right)^k
    = \frac{k!}{2\pi i A(N)} \int_{c-i\infty}^{c+i\infty}
    \frac{N^s}{s} \mathcal{Q}_k(s) \dd s,
\end{align}
where $\mathcal{Q}_k(s)$ is the coefficient of $z^k$ in the Taylor
expansion (in $z$) of
\[
    \frac1{(1+\mu z/s)(1-e^z\mathcal{P}(s+\mu z))}.
\]

While all steps can be easily justified (as done below), we cannot
directly apply Delange's Tauberian theorem to $\mathcal{Q}_k(s)$
here since each $\mathcal{Q}_k$ (except $\mathcal{Q}_0(s))$ is not a
proper Dirichlet series, but involves additional powers of $s^{-1}$.
This can be resolved as follows.

\paragraph{Shifting-the-mean at the coefficients level.} We look at
the ``translation" of the preceding parameter-shift at the
coefficient level. By definition
\begin{align*}
    A(N)\mathbb{E}\left(e^{(Y_N-\mu\log N)z}\right)
    &= \sum_{n\le N} \sum_{m\ge0} a_m(n) e^{(m-\mu\log N)z}\\
    &= \sum_{n\le N} \sum_{m\ge0} a_m(n) e^{(m-\mu\log n)z-\mu z
    \log(N/n)}.
\end{align*}
Let
\[
    b_k(n) := \sum_{m\ge0} a_m(n) (m-\mu\log n)^k.
\]
Then, by taking the coefficients of $z^k$ on both sides, we obtain
\begin{align}\label{T:k-mm}
    A(N)\mathbb{E}\left(Y_N-\mu\log N\right)^k
    = \sum_{0\le \ell \le k} \binom{k}{\ell}(-\mu)^{k-\ell}
    \sum_{n\le N}b_\ell(n)\left(\log\frac Nn\right)^{k-\ell}.
\end{align}
We will see that the growth order of $\sum_{n\le N} b_k(n)$ is the
power $N^\rho$ times an additional logarithmic term; it then follows
that the weighted sum on the right-hand side is of the same order by
a simple partial summation (see below for more details).

Now observe that (assuming again that $N$ is not an integer)
\[
    \frac1{2\pi i} \int_{c-i\infty}^{c+i\infty}
    \frac{N^s}{s^m} \,\sum_{j\ge1}\alpha(j)j^{-s} \dd s
    = \frac1{(m-1)!}\sum_{n\le N} \alpha(n) \left(\log
    \frac Nn\right)^{m-1}\qquad(m=1,2,\dots),
\]
where $c$ is taken to be any real number greater than the abscissa
of absolute convergence of the function defined by the series
$\sum_{j\ge1}\alpha(j)j^{-s}$. So, this, together with
(\ref{T:k-mm}), justifies (\ref{T:fsz}).

\paragraph{A probabilistic interpretation.} Given $N$, consider a
random factorization of a number $n\le N$, namely, a random product
$p_1\cdots p_m\le n$ with all $p_j\in\mathscr{P}$ (uniformly
distributed over all $A(N)$ possible factorizations). Let $Y_N$ be
the number of factors ($=m$) and $\nu_N$ be their product ($=n$).
Then
\[
    A(N) \mathbb{E}(Y_N-\mu\log\nu_N)^k
    =\sum_{n\le N} b_k(n),
\]
which gives a probabilistic interpretation of the partial sum.

\paragraph{The Dirichlet series of $b_k(n)$.}
Define the Dirichlet series
\begin{equation}
  \label{mks}
    \mathcal{M}_k(s) := \sum_{n\ge1} b_k(n) n^{-s}
    =\sum_{n\ge1} n^{-s}\sum_{m\ge0} a_n(m) (m-\mu\log n)^k.
\end{equation}
Note that $a_m(n)>0$ implies that $m\le \log_2n$, so that
\[
    (m-\mu\log n)^k = O\left((\log n)^k\right),
\]
for all non-zero terms. Hence $\mathcal{M}_k(s)$ is absolutely
convergent when $\Re(s)>\rho$ because
$\mathcal{A}(s)=1/(1-\mathcal{P}(s))$ is. Now if $\Re(s)>\rho$ and
$|z|$ is sufficiently small, then, by \eqref{T:bgf},
\begin{align}\label{T:Mks}
    \sum_{k\ge0} \frac{\mathcal{M}_k(s)}{k!} z^k
    = \sum_{n\ge1} n^{-s}\sum_m a_{n}(m) e^{(m-\mu \log n)z}
    = \frac1{1-e^z\mathcal{P}(s+\mu z)}.
\end{align}
With these $\mathcal{M}_k(s)$, the generating function
$\mathcal{Q}_k(s)$ can be decomposed as
\begin{align} \label{T:Qks}
    \mathcal{Q}_k(s) = \sum_{0\le \ell \le k} \left(\frac{-\mu}{s}
    \right)^{k-\ell} \frac{\mathcal{M}_\ell(s)}{\ell!}.
\end{align}

Our strategy will then to apply Delange's Tauberian theorem to
$\mathcal{M}_k$ for even $k$ and some auxiliary Dirichlet series for
odd $k$, and then the asymptotics of the $k$-th central moment can
be obtained easily since terms with index $\ell<k$ in (\ref{T:Qks})
will be asymptotically negligible. Indeed, we will use
(\ref{T:k-mm}).

\paragraph{Recurrence of $\mathcal{M}_k(s)$.} We now focus on
properties of $\mathcal{M}_k(s)$. By writing (\ref{T:Mks}) in the
form
\[
    \left(1-e^z\mathcal{P}(s+\mu z)\right) \sum_{\ell\ge0}
    \frac{\mathcal{M}_\ell(s)}{\ell!}\,z^\ell = 1,
\]
we see that $\mathcal{M}_k(s)$ satisfies the recurrence
\begin{align}\label{T:Mks-rr}
    \mathcal{M}_k(s) = \frac1{1-\mathcal{P}(s)} \sum_{0\le j<k}
    \binom{k}{j} \mathcal{M}_j(s) \mathcal{B}_{k-j}(s) \qquad(k\ge1),
\end{align}
with $\mathcal{M}_0(s)=1/(1-\mathcal{P}(s))$, where
\[
    \mathcal{B}_k(s) := \sum_{0\le \ell \le k}
    \binom{k}{\ell} \mu^\ell \mathcal{P}^{(\ell)}(s).
\]
For example, $\mathcal{M}_1(s)=\mathcal{B}_1(s)/(1-\mathcal{P}(s))^2
=(\mathcal{P}(s)+\mu\mathcal{P}'(s))/(1-\mathcal{P}(s))^2$.

Note that each $\mathcal{B}_k(s)$ is analytic for $\Re(s)>\kappa$
and, in particular, for $\Re(s)\ge \rho$. Moreover, the crucial
property here is
\[
    \mathcal{B}_1(\rho)= \mathcal{P}(\rho)+\mu
    \mathcal{P}'(\rho) = 1-1=0,
\]
by our construction. Similarly,
\[
    \mathcal{B}_2(\rho)= \mathcal{P}(\rho)+2\mu
    \mathcal{P}'(\rho)+\mu^2\mathcal{P}''(\rho)
    = \mu^2 \mathcal{P}''(\rho)-1 = \sigma^2/\mu,
\]
and $\mathcal{B}_1'(\rho)=\sigma^2/\mu^2$.

On the other hand, by (\ref{T:Mks-rr}), we see that
$\mathcal{M}_k(s)$ is analytic for $\Re(s)>\rho$ and for
$\Re(s)=\rho$ except at $s=\rho$. Furthermore, because
$\mathcal{B}_1(\rho)=0$, it follows by induction from
(\ref{T:Mks-rr}), that at $s=\rho$, $\mathcal{M}_k(s)$ has a pole of
order at most $\lfloor k/2\rfloor+1$.

\paragraph{Even moments.} More precisely, for even $k=2\ell$, we get
by induction \
\[
    \mathcal{M}_k(s) \sim c_k (s-\rho)^{-k/2-1} ,
\]
where
\[
    c_k = \binom{k}{2} \mu \mathcal{B}_2(\rho)c_{k-2}
    =\frac{k(k-1)}{2} \sigma^2 c_{k-2},
\]
with $c_0=\mu$, which is solved to be
\[
    c_k = \mu \left(\frac{\sigma^2}{2}\right)^{k/2} k!.
\]
We now apply Delange's Tauberian theorem and obtain
\begin{align}
    \mathbb{E}(Y_N-\mu\log\nu_N)^k
    &=\frac1{A(N)}\sum_{n\le N} b_k(n) \nonumber\\
    &\sim \frac{c_k}{\rho\Gamma(k/2+1)A(N)}
    \,N^\rho(\log N)^{k/2}\nonumber \\
    &\sim \frac{k!}{2^{k/2}(k/2)!}\sigma^k(\log N)^{k/2}.
    \label{T:nuN-even}
\end{align}

\paragraph{Odd moments.} Let now $k=2\ell-1$, $\ell\ge1$. Since the
coefficients $b_k(n)$ are not necessarily nonnegative, we cannot
directly apply Delange's Tauberian theorem. Instead, we consider the
following two auxiliary Dirichlet series
\[
    D_1(s) := \sum_{n\ge1} n^{-s} \sum_{m\ge0} a_m(n)\left(
    (m-\mu\log n)^k + (\log n)^{k/2}\right)^2,
\]
and, see \eqref{mks} and \eqref{aks},
\begin{align*}
    D_2(s) &:= \sum_{n\ge1} n^{-s} \sum_{m\ge0} a_m(n)
    \left((m-\mu\log n)^{2k} + (\log n)^{k}\right)\\
    &= \mathcal{M}_{2k}(s) +(-1)^k \mathcal{A}^{(k)}(s).
\end{align*}
The two Dirichlet series have only nonnegative coefficients, and we
will show that Delange's Tauberian theorem can be applied to both
series. The leading terms will cancel and we will have
\begin{align} \label{T:tb1}
    \frac1{A(N)}\sum_{n\le N} b_k(n)(\log n)^{k/2}
    = o\left((\log N)^k\right).
\end{align}
From this, we use the monotonicity of $(\log n)^{k/2}$ and
elementary arguments to recover the desired estimate
\begin{align} \label{T:tb2}
    \mathbb{E}(Y_N-\mu\log\nu_N)^k
=  \frac1{A(N)} \sum_{n\le N} b_k(n) = o\left((\log N)^{k/2}\right).
\end{align}

\paragraph{Proof of (\ref{T:tb1}).} Let
\[
    D_3(s) := \sum_{n\ge1} (\log n)^{k/2}b_k(n)n^{-s} .
\]
Then $D_1(s) = D_2(s) + 2D_3(s)$. By the discussions above, we can
apply Delange's theorem to $D_2(s)$ and obtain
\begin{align} \label{T:D2s}
    \frac1{A(N)}\sum_{n\le N}\sum_{m\ge0} a_m(n)
    \left((m-\mu\log n)^{2k} + (\log n)^{k}\right)
    \sim C_k (\log N)^k,
\end{align}
where $C_k = (2k)!\sigma^{2k}/(2^kk!)+1$ (this value is however
immaterial).

We now show that the partial sum of the coefficients of $D_1(s)$ has
asymptotically the same dominant term. We start from the
representation ($k=2\ell-1$)
\[
    D_3(s) = (-1)^\ell \pi^{-1/2} \int_0^\infty
    \mathcal{M}_k^{(\ell)}(s+t) t^{-1/2} \dd t,
\]
for $\Re s>\rho$, because $(-1)^\ell \mathcal{M}_k^{(\ell)}(s) =
\sum_{n\ge1} b_k(n) (\log n)^{\ell} n^{-s}$ and
\begin{align*}
    (-1)^\ell \int_0^\infty
    \mathcal{M}_k^{(\ell)}(s+t) t^{-1/2} \dd t
    &= \sum_{n\ge1} b_k(n) (\log n)^\ell n^{-s}
    \int_0^\infty e^{-t\log n} t^{-1/2} \dd t \\
    &= \Gamma(\tfrac12) \sum_{n\ge1} b_k(n) (\log n)^{k/2} n^{-s}\\
    &= \sqrt{\pi} D_3(s).
\end{align*}
We now consider the local behavior of $D_3(s)$ when $s\sim\rho$.
First, $\mathcal{M}_k(s)$ has a pole at $s=\rho$ with leading term
$c_k'(s-\rho)^{-(k+1)/2}$, for some $c_k'$. Thus
$\mathcal{M}_k^{(\ell)}(s)$ has a pole with local behavior
$c_k''(s-\rho)^{-k-1}$. It follows that for small $|w|$ and
$\Re(w)>0$,
\begin{align*}
    D_3(\rho+w) &= (-1)^\ell\pi^{-1/2} \int_0^\infty
    \mathcal{M}_k^{(\ell)}(\rho+w+t) t^{-1/2} \dd t \\
    &= O\left(\int_0^\infty |w+t|^{-k-1} t^{-1/2}\dd t\right) \\
    &= O\left(\int_0^{|w|} |w|^{-k-1} t^{-1/2} \dd t
    + \int_{|w|}^\infty t^{-k-3/2} \dd t \right)\\
    &= O\left(|w|^{-k-1/2}\right) \\
    &= o\left(|w|^{-k-1}\right).
\end{align*}
Since $D_1(s)=D_2(s)+D_3(s)$ has all coefficients nonnegative, we
can now apply Delange's theorem to $D_1(s)$ and conclude that
\[
    \frac1{A(N)}\sum_{n\le N} \sum_{m\ge0} a_m(n)\left(
    (m-\mu\log n)^k + (\log n)^{k/2}\right)^2
    \sim C_k (\log N)^k.
\]
This, together with (\ref{T:D2s}), proves (\ref{T:tb1}).

\paragraph{Proof of (\ref{T:tb2}).} Let
\[
    B_k(x) := \sum_{n\le x} b_k(n) (\log n)^{k/2},
\]
and for $x\ge2$, $\varphi(x) := (\log x)^{-k/2}$. Then
\begin{align*}
    \int_2^N B_k(x) \varphi'(x) \dd x
    &= \sum_{2\le n\le N} b_k(n) (\log n)^{k/2}
    \int_n^N \varphi'(x) \dd x \\
    &= \sum_{2\le n\le N} b_k(n) (\log n)^{k/2}
    \left((\log N)^{-k/2} -(\log n)^{-k/2}\right) \\
    &= B_k(N)(\log N)^{-k/2} - \sum_{2\le n\le N} b_k(n).
\end{align*}
Thus, by (\ref{T:tb1}),
\begin{align*}
    \sum_{1\le n\le N} b_k(n) &= b_k(1) + B_k(N)(\log N)^{-k/2}
    -\int_2^N B_k(x) \varphi'(x) \dd x\\
    &= O(1) + o\left(N^\rho (\log N)^{k/2}\right)
    + \frac k2\int_2^N B_k(x) x^{-1} (\log x)^{-k/2-1} \dd x\\
    &=O(1) + o\left(N^\rho (\log N)^{k/2}\right)
    +o\left(\int_2^N x^{\rho-1}(\log x)^{k/2-1}\dd x \right)\\
    &=o\left(N^\rho(\log N)^{k/2}\right),
\end{align*}
as required.

\paragraph{From $Y_N-\mu\log \nu_N$ to $Y_N-\mu\log N$.} The two
estimates (\ref{T:nuN-even}) and (\ref{T:tb2}) imply
\begin{align} \label{T:nuN-cm}
    \mathbb{E}\left(\frac{Y_N-\mu\log\nu_N}
    {\sigma\sqrt{\log N}}\right)^k\to
    \left\{\begin{array}{ll}
        \displaystyle \frac{k!}{(k/2)!2^{k/2}},
        &\text{if $k$ is even}\\
        0,& \text{if $k$ is odd},
    \end{array}\right.
\end{align}
which in turn implies, by the method of moments, that
\[
    \frac{Y_N-\mu\log\nu_N}{\sigma\sqrt{\log N}}
    \stackrel{d}{\to} \mathscr{N}(0,1).
\]

Our final task is to prove the same asymptotics (\ref{T:cm}) from
the two estimates (\ref{T:nuN-even}) and (\ref{T:tb2}). To that
purpose, define $S_k(x)=0$ if $x<2$ and
\[
    S_k(x) := \sum_{n\le x} b_k(n)\qquad(x\ge2).
\]
We use (\ref{T:k-mm}) and the cruder estimates (by
(\ref{T:nuN-even}), (\ref{T:tb2}) and \eqref{T:AN})
\[
    S_\ell(x) = O\left(A(x) (\log x)^{\ell/2}\right)
    = O\left(x^\rho (\log x)^{\ell/2}\right),
\]
for $\ell=0,\dots,k-1$. With this, we have
\begin{align*}
    &\sum_{0\le \ell<k}\binom{k}\ell (-\mu)^{k-\ell}
    \sum_{n\le N} b_\ell(n) \left(\log\frac Nn\right)^{k-\ell}\\
    &\qquad =O\left((\log N)^{k}+
    \sum_{0\le \ell<k}\binom{k}\ell \mu^{k-\ell}
    \int_2^N\left(\log\frac Nx\right)^{k-\ell} \dd S_\ell(x)
    \right).
\end{align*}
Now for each $\ell=0,\dots,k-1$,
\begin{align*}
    \int_2^N\left(\log\frac Nx\right)^{k-\ell} \dd S_\ell(x)
    &= O\left(\int_2^N\left(\log\frac Nx\right)^{k-\ell-1} x^{-1}
    S_\ell(x) \dd x \right) \\
    &= O\left(\int_2^N\left(\log\frac Nx\right)^{k-\ell-1}
    x^{\rho-1}(\log x)^{\ell/2} \dd x \right).
\end{align*}
Splitting the integral at $x=N/2$, and making the change of
variables $x\mapsto N/x$ for the first half, we see that
\begin{align*}
    &\int_2^N\left(\log\frac Nx\right)^{k-\ell-1}
    x^{\rho-1}(\log x)^{\ell/2} \dd x \\
    &\qquad= O\left(N^\rho \int_2^{N/2}x^{-\rho-1}
    (\log x)^{k-1-\ell}\left(\log\frac Nx\right)^{\ell/2}\dd x +
    \int_{N/2}^N x^{\rho-1} (\log x)^{\ell/2} \dd x \right)\\
    &\qquad= O\left(N^\rho (\log N)^{\ell/2}\right)\\
    &\qquad= o\left(N^\rho (\log N)^{k/2}\right),
\end{align*}
for $0\le \ell\le k-1$. This proves that
\[
     \mathbb{E}\left(Y_N-\mu\log N\right)^k
     = \frac1{A(N)} \sum_{n\le N} b_k(n) +
     o\left((\log N)^{k/2}\right),
\]
and thus the estimates in (\ref{T:cm}) hold by (\ref{T:nuN-even})
and (\ref{T:tb2}).

\paragraph{An alternative argument.} To bridge (\ref{T:nuN-cm}) and
(\ref{T:cm}), we can also argue as follows. Consider the sum
\[
    \mathbb{E}\left(\log N-\log\nu_N\right)^k
    = \frac1{A(N)}\sum_{n\le N} a(n) \left(\log\frac Nn\right)^k,
\]
which is $O(1)$ by a similar summation by parts argument as used
above. Then for even $k$
\[
    ||\log N-\log\nu_N||_k = O(1).
\]
By H\"older's inequality, this holds true also for every $k\ge0$.
Consequently, using again H\"older's inequality, we deduce
(\ref{T:cm}).

\section{Conclusions and additional remarks}

While a direct application of Tauberian theorems leads to results of
the form
\[
    \mathbb{E}(Y_N) \sim \mu \log N,
\]
we indeed prove, still relying on Tauberian arguments, that
\[
    \mathbb{E}(Y_N) = \mu \log N + o(\sqrt{\log N}),
\]
(a special case of Proposition~\ref{T:pp1}). This shows the power of
our approach. However, the estimates (\ref{T:cm}) we derived are not
strong enough so as to prove more effective bounds such as the
convergence rate to normality (or the Berry-Esseen bound).

Another corollary to our moment convergence result is the following
asymptotic approximations to all absolute central moments
\begin{align*}
    \mathbb{E}|Y_N-\mu\log N|^\beta
    &\sim 2^{\beta/2}\pi^{-1/2} \Gamma\left(\frac{\beta+1}2
    \right)(\log N)^{\beta/2},
\end{align*}
for all $\beta\ge 0$, which seem difficult to get directly from
Dirichlet series.

On the other hand, when $z\in (-\log \mathcal P(\kappa),\infty)$,
one can apply directly Delange's Tauberian theorem to the generating
function
\[
    \frac1{1-e^z P(s)},
\]
(instead of to the Dirichlet series of higher moments obtained above
by Taylor expansions in $z$); this results in the asymptotic
approximation
\[
    \mathbb{E}\left(e^{zY_N}\right)
    \sim \frac{\rho P'(\rho)}
    {\rho(z)e^z P'(\rho(z))}\,N^{\rho(z)-\rho},
\]
where $\rho(z)$ solves the equation $1=e^zP(\rho(z))$ with
$\rho(0)=\rho$. From this approximation, one might expect asymptotic
normality by straightforward argument. However, the asymptotic
result so obtained holds only pointwise, and the uniformity in $z$
is missing here. While the gap of uniformity may perhaps be filled
by applying suitable Tauberian theorems with remainders, the use of
Delange's Tauberian theorems is computationally simpler and
technically less involved.

It is clear from our proof that Theorem~\ref{T:thm} actually holds
for any Dirichlet series $\mathcal{P}(s)$ with nonnegative
coefficients and satisfying the conditions of Theorem~\ref{T:thm}.
Thus the restriction of $\mathscr{P}$ to a subset of positive
integers is not essential. For example, one can consider the
\emph{ordered totient factorizations} with $\mathcal{P}(s) =
\sum_{n\ge3} \phi(n)^{-s}$, where $\phi(n)$ is Euler's totient
function, namely, the number of positive integers $\le n$ and
relatively prime to $n$. In this case, $\kappa=1$ and $\rho\approx
2.26386$ since
\[
    \sum_{n\ge1} \phi(n)^{-s} = \zeta(s) \prod_{p\,:\,
    \text{prime}} \left(1-p^{-s} + (p-1)^{-s}\right);
\]
see \cite{BM90} for a detailed studied of this Dirichlet series.

How to compute $\rho$ to high degree of precision? In general, the
problem is not easy, for example, if $\mathcal{P}(s) = \sum_{n\ge2}
\lceil n^\beta(\log n)^c\rceil^{-s}$ for $\beta, c>0$; the case of
totient factorization is similar. The easy cases are when
$\mathcal{P}(s)$ can be expressed in terms of $\zeta$-functions such
as $\mathcal{P}(s)=\zeta(s)-1$ (all integers $>1$) or
$\mathcal{P}(s)= \zeta(s)/\zeta(2s)-1$ (square-free integers $>1$).
Take now $\mathcal{P}(s)= \sum_{p\,:\,\text{prime}}p^{-s}$. The zero
of $\mathcal{P}(s)=1$ can also be easily computed by using the
relation (see \cite{FV96})
\[
    \sum_{p\,:\,\text{prime}}p^{-s} = \sum_{k\ge1}
    \frac{\mu(k)}{k}\, \log\zeta(ks)\qquad(\Re(s)>1),
\]
where $\mu(k)$ is M\"obius function. This readily gives
\[
    \rho \approx 1.39943\,33287\,26330\,31820\,28072\dots,
\]
and
\begin{align*}
    \mu  &\approx 0.57764\,86251\,95138\,05440\,61351\dots,\\
    \sigma^2 &\approx 0.48439\,65045\,13598\,28128\,07456\dots.
\end{align*}

We indicated a few directions to which our approach can be extended
in Introduction. But can a similar idea be modified so as to deal
with arithmetic functions with mean other than $\log N$?

\end{document}